\documentclass[article]{siamart220329}

\usepackage{amsfonts}
\usepackage{amsmath}
\usepackage{graphicx}
\usepackage{subfig}
\usepackage{epstopdf}
\usepackage{hyperref}
\ifpdf%
  \DeclareGraphicsExtensions{.eps,.pdf,.png,.jpg}
\else
  \DeclareGraphicsExtensions{.eps}
\fi
\usepackage{amsopn}
\DeclareMathOperator{\diag}{diag}
\usepackage{booktabs}

\newcommand{\TheTitle}{Local Fourier Analysis of P-Multigrid for High-Order Finite Element Operators}

\newcommand{\TheShortTitle}{LFA of P-Multigrid for High-Order FEM}

\newcommand{\TheFunding}{This work is supported by the Exascale Computing Project (17-SC-20-SC), a collaborative effort of two U.S. Department of Energy organizations (Office of Science and the National Nuclear Security Administration) responsible for the planning and preparation of a capable exascale ecosystem, including software, applications, hardware, advanced system engineering and early testbed platforms, in support of the nation’s exascale computing imperative.}

\author{
  Jeremy L Thompson\thanks{Department of Computer Science, University of Colorado, Boulder, CO
    (\email{jeremy@jeremylt.org}).}
  \and Jed Brown\thanks{Department of Computer Science, University of Colorado, Boulder, CO
    (\email{jed@jedbrown.org}).}
  \and Yunhui He\thanks{Department of Applied Mathematics, University of Waterloo, Waterloo, ON, Canada
    (\email{yunhui.he@uwaterloo.ca}).}
}
\newcommand{\TheShortAuthors}{
  J. L. Thompson, J. Brown, and Y. He
}

\title{{\TheTitle}\thanks{\TheFunding}}
\headers{\TheShortTitle}{\TheShortAuthors}

\ifpdf
\hypersetup{
  pdftitle={\TheTitle},
  pdfauthor={\TheShortAuthors}
}
\fi

\begin{document}

\maketitle

\vspace{1cm}

\begin{abstract}
Multigrid methods are popular for solving linear systems derived from discretizing PDEs.
Local Fourier Analysis (LFA) is a technique for investigating and tuning multigrid methods.
P-multigrid is popular for high-order or spectral finite element methods, especially on unstructured meshes.
In this paper, we introduce LFAToolkit.jl, a new Julia package for LFA of high-order finite element methods.
LFAToolkit.jl analyzes preconditioning techniques for arbitrary systems of second order PDEs and supports mixed finite element methods.
Specifically, we develop LFA of p-multigrid with arbitrary second-order PDEs using high-order finite element discretizations and examine the performance of Jacobi and Chebyshev smoothing for two-grid schemes with aggressive p-coarsening.
A natural extension of this LFA framework is the analysis of h-multigrid for finite element discretizations or finite difference discretizations that can be represented in the language of finite elements.
With this extension, we can replicate previous work on the LFA of h-multigrid for arbitrary order discretizations using a convenient and extensible abstraction.
Examples in one, two, and three dimensions are presented to validate our LFA of p-multigrid for the Laplacian and linear elasticity.
\end{abstract}

\begin{keywords}
  local Fourier analysis, p-multigrid, high-order, finite elements
\end{keywords}

\section{Introduction}\label{sec:intro}

Multigrid methods \cite{brandt1982guide, briggs2000multigrid, stuben1982multigrid} are popular for solving linear systems derived from discretizing PDEs.
Local Fourier Analysis (LFA) \cite{brandt1977multi, wienands2004practical} is a powerful tool for predicting multigrid performance and tuning of multigrid components by examining the spectral radius and/or condition number of the {\em{symbol}} of the underlying operator, which enables sharp predictions for large-scale problems.

High-order finite element methods offer accuracy advantages over low-order finite elements \cite{demkowicz1989toward, oden1989toward, rachowicz1989toward}; however, high-order finite elements are less common because a linear operator or the Jacobian of a non-linear operator rapidly loses sparsity in a sparse matrix representation.
Matrix-free formulations of these methods \cite{brown2010efficient, deville2002highorder, knoll2004jacobian}, provide efficient implementations of these methods on modern hardware \cite{libceed-joss, fischer2020scalability, kronbichler2019multigrid}.

LFA of h-multigrid with high-order finite elements was developed in \cite{he2020two} for Lagrange bases with uniformly spaced nodes.
P-multigrid, developed by R{\o}nquist and Patera \cite{ronquist1987spectral}, is based on decreasing the order of the bases in high-order or spectral finite element methods rather than aggregating elements to coarsen the mesh.
While it is possible \cite{davydov2019matrix} to use high-order finite element methods with matrix-free implementation of h-multigrid preconditioning on complex problems, p-multigrid can be a more natural fit for high-order methods on unstructured meshes.
The benefits of using $p$-multigrid for matrix-free operators on unstructured meshes for solid mechanics engineering applications is discussed in \cite{ratel-preprint}.

In this paper, we develop LFA of p-multigrid with arbitrary second-order PDEs using high-order finite element discretizations and investigate Jacobi and Chebyshev smoothing for two-grid schemes by way of LFAToolkit.jl \cite{thompson2021toolkit}, a new Julia package for LFA of high-order finite element methods.
Several software packages have been developed for LFA of h-multigrid methods \cite{rittich2018extending,kahl2020automated,wienands2004practical}, however, to the best of our knowledge, no packages implement LFA of p-multigrid methods, especially with arbitrary PDEs.
Although \cite{van2011discrete} discussed LFA of p-multigrid for the discontinuous Galerkin method, this formulation cannot be extended to the continuous Galerkin method.
Our LFA formulation for p-multigrid with high-order finite element discretizations for the continuous Galerkin method can be generalized to reproduce previous work on h-multigrid for high-order finite elements \cite{he2020two} and can be extended to LFA for h-multigrid of finite difference discretizations that can be represented via finite element discretizations as well as LFA of block smoothers \cite{maclachlan2011overlapping,brown2019local}.
Our work here could be represented in terms of the aLFA framework \cite{kahl2020automated}; however we choose a representation that will be familiar for the high-order finite element and spectral elements community.
We give details on the mathematical representation of LFA of grid-transfer operators for p-multigrid for the continuous Galerkin method, which seems absent in current literature.
Moreover, our LFA framework for p-multigrid generalizes and unifies classical LFA for h-multigrid, which would be a rich area for research.

We investigate the performance of Jacobi and Chebyshev semi-iterative smoothers for p-multigrid with aggressive coarsening for the scalar Laplacian in one and two dimensions, and we validate our LFA of p-multigrid with the scalar Laplacian in three dimensions against numerical experiments.
Our analysis demonstrates that the performance of p-multigrid with these two smoothers degrades as we coarsen more aggressively.
Finally, we analyze the performance of p-multigrid with a Chebyshev smoother for three dimensional linear elasticity to demonstrate ability of our LFA framework to handle more complex PDEs.

This paper is organized as follows.
In \cref{sec:notation} we outline the notation for LFA used in this paper.
In \cref{sec:lfa} we develop the notation for LFA of second-order PDEs with arbitrary order bases, dimension, and number of components used in LFAToolkit.jl, and we use this notation to develop LFA of p-multigrid and the Jacobi and Chebyshev smoothers.
\Cref{sec:results} contains numerical results investigating the performance of Jacobi and Chebyshev smoothing for p-multigrid on the one, two, and three dimensional scalar Laplacian and three dimensional linear elasticity.
Finally, \cref{sec:conclusion} contains concluding remarks.

\subsection{Reproducibility}\label{sec:reproducibility}

The numerical results for LFA in this paper were generated using the Julia package LFAToolkit.jl \cite{thompson2021toolkit}.
This package is under active development and may be found on GitHub at \href{https://github.com/jeremylt/LFAToolkit.jl}{jeremylt/LFAToolkit.jl}.
This repository contains Julia scripts and interactive Jupyter notebooks that can replicate the tables and plots in this paper.

Our numerical experiments demonstrating actual convergence rates of p-multigrid were conducted using libCEED \cite{libceed-joss} with PETSc \cite{petsc-user-ref} p-multigrid example found in the libCEED repository.

\section{Definitions and Notation}\label{sec:notation}

Consider a scalar Toeplitz operator $L_h$ on an infinite one dimensional uniform nodal grid $G_h$,
\begin{equation}
\begin{split}
L_h \mathrel{\hat{=}} \left[ s_\kappa \right]_h \left( \kappa \in V \right);\\
L_h w_h \left( x \right) = \sum_{\kappa \in V} s_\kappa w_h \left( x + \kappa h \right),
\end{split}
\end{equation}
where $V \subset \mathbb{Z}$ is a finite index set, $s_\kappa \in \mathbb{R}$ are constant coefficients, and $w_h \left( x \right)$ is a $l^2$ function on $G_h$. In the terminology of stencils, $s_{\kappa}$ are stencil weights that are nonzero on the neighborhood $\kappa \in V$.

Since $L_h$ is Toeplitz, it can be diagonalized by the standard Fourier modes $\varphi \left( \theta, x \right) = e^{\imath \theta x / h}$, with $i^2 = -1$, which alias according to $\varphi(\theta+ 2\pi, x) = \varphi(\theta, x)$ on all grid points $x \in h \mathbb Z$.

\begin{definition}[Symbol of $L_h$]\label{def:symbol}
If for all grid functions $\varphi \left( \theta, x \right)$ we have
\begin{equation}
L_h \varphi \left( \theta, x \right) = \tilde{L}_h \left( \theta \right) \varphi \left( \theta, x \right),
\end{equation}
then we define $\tilde{L}_h \left( \theta \right) = \sum_{\kappa \in V} s_\kappa e^{\imath \theta \kappa}$ as the symbol of $L_h$, where $\imath^2 = -1$.
\end{definition}

This definition can be extended to a $q \times q$ linear system of operators by
\begin{equation}
\mathbf{L}_h =
\begin{bmatrix}
    L_h^{1, 1} && \cdots && L_h^{1, q}        \\
    \vdots               && \vdots && \vdots  \\
    L_h^{q, 1} && \cdots && L_h^{q, q}        \\
\end{bmatrix},
\end{equation}
where $L_h^{i, j}$, $i, j \in \lbrace 1, 2, \dots, q \rbrace$ are given by scalar Toeplitz operators describing how component $j$ appears in the equation for component $i$.
The components here may represent different fields and/or nodes of a finite element basis.
The symbol of $\mathbf{L}_h$, denoted $\tilde{\mathbf{L}}_h$, is a $q \times q$ matrix-valued function of $\theta$ given by $\left( \tilde{\mathbf{L}}_h \right)_{i, j} = \tilde{L}_h^{i, j} \left( \theta \right)$.
For a system of equations representing an error propagation operator in a relaxation scheme, the spectral radius of the symbol matrix determines how rapidly the scheme decreases error at a target frequency.

These definitions are readily extended to $d$ dimensions by taking the neighborhood $V \subset \mathbb Z^d$ and letting $\theta \in [-\pi,\pi)^d$.
For standard coarsening from a fine grid with a mesh size $h$ to a coarse grid with a mesh size of $H = 2 h$, low frequencies are given by $\theta \in T^{\text{low}} = \left[ - \pi / 2, \pi / 2 \right)^d$ and high frequencies are given by $\theta \in \left[-\pi, \pi \right)^d \setminus T^{\text{low}}$ or equivalently by periodicity, $\theta \in T^{\text{high}} = \left[ - \pi / 2, 3 \pi / 2 \right)^d \setminus T^{\text{low}}$.

\section{Local Fourier Analysis for P-Multigrid}\label{sec:lfa}

We now develop the LFA formulation used in LFAToolkit.jl, first in one dimension and then in multiple dimensions, followed by LFA of polynomial smoothers and LFA of p-multigrid with high-order finite elements.

\subsection{High-Order Finite Elements}\label{sec:highorder}

We will use the representation of the weak form of linear second-order PDEs described in \cite{brown2010efficient}, which is given by
\begin{equation}
\langle v, f \left( u \right) \rangle = \int_{\Omega}
\begin{bmatrix}
  v^T & \nabla v^T    \\
\end{bmatrix}
\begin{bmatrix}
  f_{0, 0} & f_{0, 1} \\
  f_{1, 0} & f_{1, 1} \\
\end{bmatrix}
\begin{bmatrix}
  u                   \\
  \nabla u            \\
\end{bmatrix}
= \int_{\Omega} f v,\quad \forall v \in \mathcal V
\end{equation}
for some suitable $\mathcal V \subseteq H_0^1 \left( \Omega \right)$.
In this equation, $f_{i, j}$ may come from a linear PDE or the linearization of a non-linear problem.
Boundary terms have been omitted, as they are not present on the infinite uniform grid $G_h$.

This omission of boundary conditions corresponds to periodic boundary conditions on a finite domain; however, in practice LFA on the infinite uniform grid offers good predictions for finite domains with other boundary conditions \cite{wienands2004practical}.

Selecting a finite element basis, we can discretize this weak form and produce
\begin{equation}\label{pdediscrete}
\mathbf{A} \mathbf{u} = \mathbf{b}.
\end{equation}

Using the algebraic representation of PDE operators given in \cite{brown2010efficient}, the PDE operator $\mathbf{A}$ is of the form
\begin{equation}\label{efficienthighorder}
\mathbf{A} = \mathbf{G}^T \mathbf{A}_e \mathbf{G},\quad \text{with} \,\,\mathbf{A}_e = \mathbf{B}^T \mathbf{D} \mathbf{B},
\end{equation}
where $\mathbf{G}$ represents the element assembly operator, $\mathbf{B}$ is a basis operator which computes the values and derivatives of the basis functions at the quadrature points, and $\mathbf{D}$ is a block diagonal operator which provides the pointwise application of the bilinear form on the quadrature points, to include quadrature weights and the change in coordinates between the physical and reference space.

Consider the specific case of a block Toeplitz operator representing an arbitrary second order scalar PDE in one dimension with Lagrange basis of polynomial degree $p$ given in the form of \cref{pdediscrete} and \cref{efficienthighorder}.
The basis $\mathbf{B}$ is defined on the mesh grid with points given by $x_i$, for $i \in \left\lbrace 1, 2, \dots, p + 1 \right\rbrace$ on one element.
The nodes on the left and right boundaries of the element map to the same Fourier mode when localized to nodes unique to a single finite element, so we can compute the symbol matrix of size $p \times p$ as
\begin{equation}\label{symbolscalar1d}
\tilde{\mathbf{A}} = \mathbf{Q}^T \left( \mathbf{A}_e \odot \left[ e^{\imath \left( x_j - x_i \right) \theta / h} \right] \right) \mathbf{Q},
\end{equation}
where $\odot$ represents pointwise multiplication of the elements, $h$ is the length of the element, and $i, j \in \lbrace 1, 2, \dots, p + 1 \rbrace$.
In the pointwise product $\mathbf{A}_e \odot \left[ e^{\imath \left( x_j - x_i \right) \theta / h} \right]$, the $\left( i, j \right)$ entry is given by $\left( \mathbf{A}_e \right)_{i, j} e^{\imath \left( x_j - x_i \right) \theta / h}$.
$\mathbf{Q}$ is a $\left( p + 1 \right) \times p$ matrix that localizes Fourier modes to each finite element, given by
\begin{equation}
\mathbf{Q} =
\begin{bmatrix}
    \mathbf{I}   \\
    \mathbf{e}_0 \\
\end{bmatrix} =
\begin{bmatrix}
    1      && 0      && \cdots && 0      \\
    0      && 1      && \cdots && 0      \\
    \vdots && \vdots && \vdots && \vdots \\
    0      && 0      && \cdots && 1      \\
    1      && 0      && \cdots && 0      \\
\end{bmatrix}.
\end{equation}
We refer to $\mathbf{Q}$ as the Fourier mode localization operator.

The computation of this symbol matrix extends to more complex PDE with multiple components and in higher dimensions.
Multiple components are supported by extending the $p \times p$ system of Toeplitz operators in \cref{symbolscalar1d} to a $\left( n \cdot p \right) \times \left( n \cdot p \right)$ system of operators, where $n$ is the number of components in the PDE.
The localization operator for a multi-component PDE is given by $\mathbf{Q}_n = \mathbf{I}_n \otimes \mathbf{Q}$, where $I_n$ is the identity matrix with size $n \times n$.
In general, we omit the subscript indicating the number of components for the localization operator.

The infinite uniform grid $G_h$ is extended into higher dimensions by taking the direct sum of the one dimensional grid.
Tensor products are used to extend the one dimensional bases into higher dimensions.
The basis evaluation operators in two dimensions are given by
\begin{equation}
\begin{split}
\mathbf{B}_{\text{interp2d}} = \mathbf{B}_{\text{interp1d}} \otimes \mathbf{B}_{\text{interp1d}}, \\
\mathbf{B}_{\text{grad2d}} =
\begin{bmatrix}
    \mathbf{B}_{\text{grad1d}} \otimes \mathbf{B}_{\text{interp1d}} \\
    \mathbf{B}_{\text{interp1d}} \otimes \mathbf{B}_{\text{grad1d}} \\
\end{bmatrix},
\end{split}
\end{equation}
where $\mathbf{B}_{interp*d}$ and $\mathbf{B}_{grad*d}$ are the basis interpolation and gradient operators, respectively.

In a similar fashion, the localization operator for Fourier modes in two dimensions is given by
\begin{equation}
\mathbf{Q}_{\text{2d}} = \mathbf{Q} \otimes \mathbf{Q}.
\end{equation}
As with tensor product bases, an analogous computation can be done in higher dimensions.
Again, we generally omit the subscript indicating the dimension of the Fourier mode localization operator.

\begin{definition}\label{def:high_order_symbol}
The symbol matrix of a finite element operator for an arbitrary second order PDE with any basis order, dimension, and number of components is given by
\begin{equation}\label{symbolhighorder}
\tilde{\mathbf{A}} = \mathbf{Q}^T \left( \mathbf{A}_e \odot \left[ e^{\imath \left( \mathbf{x}_j - \mathbf{x}_i \right) \cdot \boldsymbol{\theta} / \mathbf{h}} \right] \right) \mathbf{Q},
\end{equation}
where $\odot$ represents pointwise multiplication of the elements, $\mathbf{h}$ is the length of the element in each dimension, $\boldsymbol{\theta}$ is the target frequency in each dimension, $i, j \in \lbrace 1, 2, \dots, n \cdot \left( p + 1 \right)^d \rbrace$, $n$ is the number of components, $p$ is the polynomial degree of the discretization, and $d$ is the dimension of the finite element basis.
$\mathbf{A}_e$ is the finite element operator for the element and $\mathbf{Q}$ is the localization of Fourier modes on an element.
\end{definition}

Note that this LFA notation is applicable to any second-order PDE with a weak form that can be represented by \cref{efficienthighorder}.
This representation is used in LFAToolkit.jl, where the users provide the finite element basis $\mathbf{B}$, the node to mode mapping $\mathbf{Q}$, and the pointwise representation of the weak form $\mathbf{D}$, and the software can provide the LFA of the PDE with various preconditioners.
The nodes may be defined by points $\mathbf x_i$ or more generally, as dual basis functions that may be applied to the Fourier modes; however, we restrict ourselves to nodal bases in this paper.

If the pointwise representation of the weak form $\mathbf{D}$ has a tensor product structure, then the element operator $\mathbf{A}_e$ and therefore the symbol $\tilde{\mathbf{A}}$ will have a tensor product structure, as shown in \cite{he2020two}.
We omit this optimization in this analysis and LFAToolkit.jl in favor of supporting second-order PDEs that do not have a tensor product decomposition.

\subsection{Polynomial Smoothers}\label{sec:smooth}

Multigrid methods require a fine grid smoother, typically applied both before and after the coarse grid correction is computed.
In this section, we compute the symbol of the error propagation operator for two polynomial smoothers, Jacobi and Chebyshev.

The error propagation operator for a smoother is given by $\mathbf{S} = \mathbf{I} - \mathbf{M}^{-1} \mathbf{A}$, where $\mathbf{M}^{-1}$ is given by the particular smoother under investigation.
The LFA-predicted convergence factor is given by the maximum spectral radius of the symbol across all frequencies
\begin{equation}
\mu \left( \omega \right) = \max_{\boldsymbol{\theta} \in T} \rho \left( \tilde{\mathbf{S}} \left( \nu, \omega, \boldsymbol{\theta} \right) \right),
\end{equation}
where $ \rho \left( \tilde{\mathbf{S}} \left( \nu, \omega, \boldsymbol{\theta} \right)\right)$ denotes the spectral radius of the matrix symbol $\tilde{\mathbf{S}}$.
Our goal is to select tuning parameters for the polynomial smoothers to minimize the LFA-predicted convergence factor of two-grid p-multigrid cycles with agressive coarsening.

\subsubsection{Jacobi Smoother}\label{sec:jacobi}

With Jacobi smoothing, $\mathbf{M}^{-1}$ is given by a weighted inverse of the true operator diagonal, $\omega \diag \left( \mathbf{A} \right)^{-1}$.
Following the derivation from \cref{sec:highorder}, the symbol of the Jacobi error propagation operator therefore is given by
\begin{equation}
\tilde{\mathbf{S}} \left( \omega, \boldsymbol{\theta} \right) = \mathbf{I} - \tilde{\mathbf{M}}^{-1} \left( \omega, \boldsymbol{\theta} \right) \tilde{\mathbf{A}} \left( \boldsymbol{\theta} \right) = \mathbf{I} - \omega \left( \mathbf{Q}^T \diag \left( \mathbf{A}_e \right) \mathbf{Q} \right)^{-1} \tilde{\mathbf{A}} \left( \boldsymbol{\theta} \right),
\end{equation}
where $\tilde{\mathbf{M}}^{-1}$ has been simplified by the fact that $e^{\imath \left( \mathbf{x}_i - \mathbf{x}_i \right) \cdot \boldsymbol{\theta} / \mathbf{h}} = 1$.

If multiple pre or post-smoothing passes are used by our multigrid algorithm, we take the product of the symbol matrix with itself to represent repeated application of the Jacobi smoother.

\begin{definition}\label{def:jacobi_symbol}
The symbol of the error propagation operator for Jacobi smoothing is given by
\begin{equation}
\tilde{\mathbf{S}} \left( \nu, \omega, \boldsymbol{\theta} \right) = \left( \mathbf{I} - \omega \left( \mathbf{Q}^T \diag \left( \mathbf{A}_e \right) \mathbf{Q} \right)^{-1} \tilde{\mathbf{A}} \left( \boldsymbol{\theta} \right) \right)^\nu,
\end{equation}
where $\nu$ is the number of smoothing passes and $\omega$ is the smoothing parameter.
\end{definition}

\begin{figure}[!tbp]
  \centering
  \subfloat[Spectrum of Jacobi for $p = 2$]{\includegraphics[width=0.48\textwidth]{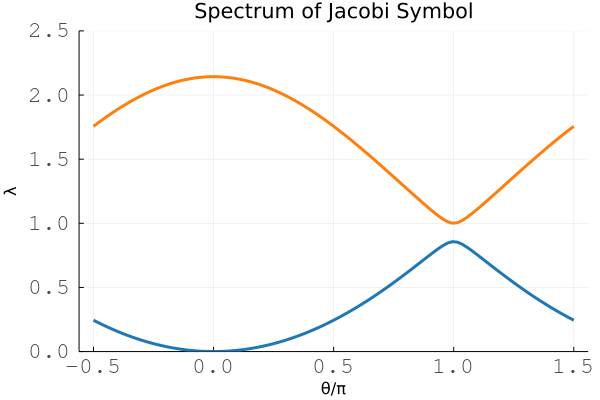}\label{fig:jacobi_spectrum_2}}
  \hfill
  \subfloat[Spectrum of Jacobi for $p = 4$]{\includegraphics[width=0.48\textwidth]{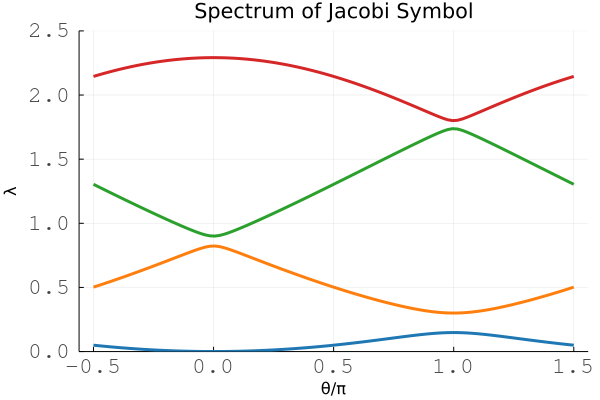}\label{fig:jacobi_spectrum_4}}
  \caption{Jacobi smoothing for high-order finite elements for the 1D Laplacian.}
\end{figure}

Using \cref{def:jacobi_symbol}, we plot the eigenvalues of $\tilde{\mathbf{M}}^{-1} \tilde{\mathbf{A}} \left( \boldsymbol{\theta} \right)$ over the interval $\theta \in \left[ - \pi / 2, 3 \pi / 2 \right)$ in \cref{fig:jacobi_spectrum_2,fig:jacobi_spectrum_4} for the one dimensional Laplacian with a 2nd and 4th degree $H^1$ Lagrange basis functions on Gauss-Lobatto points.
The eigenvalues have been colored to identify the characteristic curves.

\subsubsection{Chebyshev Smoother}\label{sec:chebyshev}

Damped Jacobi with parameter $\omega$ is based on the degree 1 monic polynomial with root at $\alpha = 1/\omega > 0$.
For $\nu=1$ smoothing passes, this is the minimal degree 1 monic polynomial on the interval $[\alpha-c, \alpha+c]$ for any positive $c$.
However, for $\nu \ge 2$, the Jacobi smoother is a degree $\nu$ monic polynomial that is not minimal on the interval (because it ``concentrates'' at the center $\alpha$).
Chebyshev smoothers are based on a stable recurrence for the polynomial of any degree that is minimial on the specified target interval of the positive real axis, which should contain the eigenvalues of all modes the smoother is responsible for.
It is well known that polynomial smoothers allow more aggressive coarsening than Jacobi \cite{brannick2015polynomial}.
For further discussion of the error propagation properties of the Chebyshev semi-iterative method, see \cite{gutknecht2002revisited}.

We use the Jacobi preconditioned operator, $\left( \diag {\mathbf{A}} \right)^{-1} {\mathbf{A}}$, in this iteration instead of the finite element operator ${\mathbf{A}}$, similar to the discussion in \cite{adams2003parallel}.

The terms in the Chebyshev semi-iterative method can be modeled by the three term recurrence relation given by
\begin{equation}
\mathbf{u}_k = - \left( \mathbf{r}_{k - 1} + \alpha \mathbf{u}_{k - 1} + \beta_{k - 2} \mathbf{u}_{k - 2} \right) / \gamma_{k - 1},
\label{eq:chebyshev_recursive}
\end{equation}
with left-preconditioned residual $\mathbf r_{k-1} = \left( \diag {\mathbf{A}} \right)^{-1} \left( \mathbf b - \mathbf A \mathbf u_{k-1} \right)$.
The target interval $[\lambda_{\min}, \lambda_{\max}] = [\alpha - c, \alpha + c]$ should be selected to be as narrow as possible while containing the eigenvalues of the preconditioned operator $\left( \diag {\mathbf{A}} \right)^{-1} \mathbf A$ that the smoother is responsible for.
(Other polynomial methods should be used for operators that are not symmetric positive definite or nearly so.)
The parameters $\beta$ and $\gamma$ are given by the recurrence
\begin{equation}
\begin{tabular}{c c}
$\beta_0 = - \frac{c^2}{2 \alpha}$ & $\gamma_0 = - \alpha$\\
$\beta_k = \frac{c}{2} \frac{T_k \left( \eta \right)}{T_{k + 1} \left( \eta \right)} = \left( \frac{c}{2} \right)^2 \frac{1}{\gamma_k}$ & $\gamma_k = \frac{c}{2} \frac{T_{k + 1} \left( \eta \right)}{T_k \left( \eta \right)} = - \left( \alpha + \beta_{k - 1} \right)$.
\end{tabular}
\end{equation}
In this equation, $T_i \left( \zeta \right) = 2 \zeta T_{i - 1} \left( \zeta \right) - T_{i - 2} \left( \zeta \right)$ are the classical Chebyshev polynomials, which are evaluated at the point $\eta = - \alpha / c$.
First order Chebyshev is equivalent to Jacobi smoothing with the classical choice of parameter value of $\omega = 1 / \alpha$, so we should expect relatively poor smoothing for $k = 1$.

The residual in the Chebyshev semi-iterative method can therefore be modeled by the three term recurrence
\begin{equation}
\mathbf{r}_k = \left( \left( \diag {\mathbf{A}} \right)^{-1} {\mathbf{A}} \mathbf{r}_{k - 1} - \alpha \mathbf{r}_{k - 1} - \beta_{k - 2} \mathbf{r}_{k - 2} \right) / \gamma_{k - 1}.
\label{eq:chebyshev_error_recursive}
\end{equation}

Using the recurrence relation given in \cref{eq:chebyshev_error_recursive}, we can define the error propagation of the $k$th order Chebyshev smoother in terms of the error in the first term:
\begin{equation}
\begin{tabular}{c}
$\mathbf{E}_0 = \mathbf{I}$\\
$\mathbf{E}_1 = \mathbf{I} - \frac{1}{\alpha} \left( \diag {\mathbf{A}} \right)^{-1} {\mathbf{A}}$\\
$\mathbf{E}_k = \left(\left( \diag {\mathbf{A}} \right)^{-1} {\mathbf{A}} \mathbf{E}_{k - 1} - \alpha \mathbf{E}_{k - 1} - \beta_{k - 2} \mathbf{E}_{k - 2} \right) / \gamma_{k - 1}$
\end{tabular}
\label{eq:chebyshev_error_propagation}
\end{equation}
With this recursive definition of the error propagation operator, we can define the symbol for Chebyshev smoothing.

\begin{definition}\label{def:chebyshev_symbol}
The symbol of the error propagation operator for a $k$th order Chebyshev smoother based on the Jacobi preconditioned operator is given by
\begin{equation}
\tilde{\mathbf{S}} \left( \nu, k, \boldsymbol{\theta} \right) = \left( \tilde{\mathbf{E}}_k \right)^\nu,
\end{equation}
where $\nu$ is the number of smoothing passes and $\tilde{\mathbf{E}}_k \left( \mathbf{\boldsymbol{\theta}} \right)$ is given by the recursive definition
\begin{equation}
\begin{tabular}{c}
$\tilde{\mathbf{E}}_0 \left( \boldsymbol{\theta} \right) = \mathbf{I}$\\
$\tilde{\mathbf{E}}_1 \left( \boldsymbol{\theta} \right) = \mathbf{I} - \frac{1}{\alpha} \tilde{\mathbf{A}}_J \tilde{\mathbf{A}} \left( \mathbf{\theta} \right)$\\
$\tilde{\mathbf{E}}_k \left( \boldsymbol{\theta} \right) = \left( \tilde{\mathbf{A}}_J \tilde{\mathbf{A}} \left( \boldsymbol{\theta} \right) \tilde{\mathbf{E}}_{k - 1} \left( \boldsymbol{\theta} \right) - \alpha \tilde{\mathbf{E}}_{k - 1} \left( \boldsymbol{\theta} \right) - \beta_{k - 2} \tilde{\mathbf{E}}_{k - 2} \left( \boldsymbol{\theta} \right) \right) / \gamma_{k - 1}$
\end{tabular}
\end{equation}
with $\tilde{\mathbf{A}}_J = \left( \mathbf{Q}^T \diag \left( \mathbf{A}_e \right) \mathbf{Q} \right)^{-1}$ giving the symbol of the Jacobi operator.
\end{definition}

\begin{figure}[!tbp]
  \centering
  \subfloat[Spectrum of Chebyshev for $p = 2$]{\includegraphics[width=0.48\textwidth]{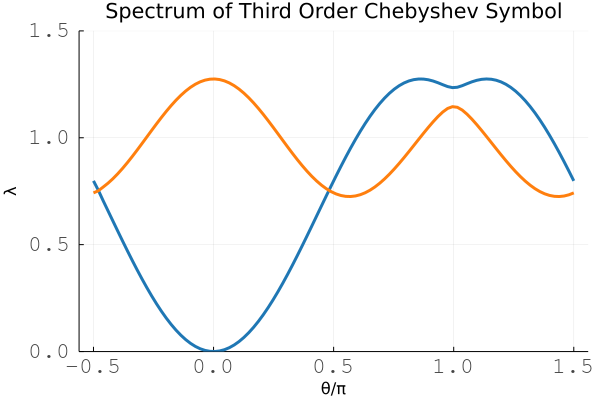}\label{fig:chebyshev_spectrum_2}}
  \hfill
  \subfloat[Spectrum of Chebyshev for $p = 4$]{\includegraphics[width=0.48\textwidth]{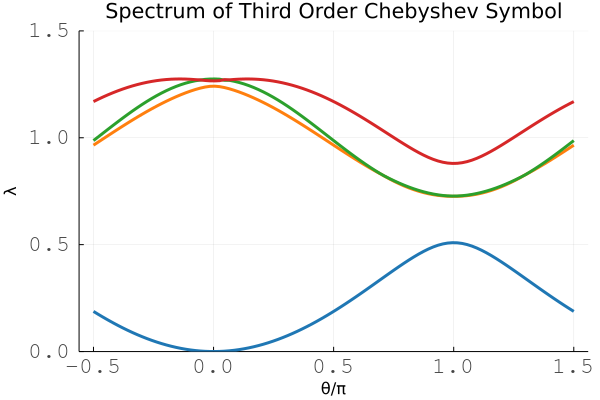}\label{fig:chebyshev_spectrum_4}}
  \caption{Chebyshev smoothing for high-order finite elements for the 1D Laplacian.}
\end{figure}

Using \cref{def:chebyshev_symbol}, we plot the eigenvalues of $\tilde{\mathbf{M}}^{-1} \tilde{\mathbf{A}} \left( \boldsymbol{\theta} \right)$ in \cref{fig:chebyshev_spectrum_2,fig:chebyshev_spectrum_4} for the one dimensional Laplacian with a 2nd and 4th degree $H^1$ Lagrange basis functions on Gauss-Lobatto points.
In both plots, we use the computation of the extremal eigenvalues discussed below.
These eigenvalues have been colored to identify the characteristic curves.
We see improved smoothing when compared to the Jacobi smoothing in \cref{fig:jacobi_spectrum_2,fig:jacobi_spectrum_4}.

The Chebyshev semi-iterative method is known to be sensitive to the quality of the estimates of the extremal eigenvalues $\lambda_{\text{min}}$ and $\lambda_{\text{max}}$.
Large eigenvalues are associated with higher frequencies for differential operators so $\lambda_{\min}$ for effective smoothing is not the minimal eigenvalue, which is hard to compute and goes to zero under grid refinement.
PETSc \cite{petsc-user-ref} estimates the eigenvalues of the preconditioned operator via Krylov iterations, yielding an estimate of the maximal eigenvalue, $\hat{\lambda}_{\text{max}}$, which is used to set the target interval according to $\lambda_{\text{min}} = 0.1 \hat{\lambda}_{\text{max}}$ and $\lambda_{\text{max}} = 1.1 \hat{\lambda}_{\text{max}}$.
This default lower bound was chosen empirically for rapid coarsening with smoothed aggregation AMG and the upper bound incorporates a safety due to $\hat{\lambda}_{\max}$ generally being an underestimate of the true maximum eigenvalue.

In LFAToolkit.jl, we also want to target the upper part of the spectrum of the error; we estimate the spectral radius of the symbol of the Jacobi preconditioned operator by sampling the the frequencies at a small number of values to estimate $\hat{\lambda}_{\text{max}}$.
We then take $\lambda_{\text{min}} = 0.1 \hat{\lambda}_{\text{max}}$ and $\lambda_{\text{max}} = 1.0 \hat{\lambda}_{\text{max}}$.

\subsection{P-Multigrid}\label{sec:multigrid}

With this representation of the symbol of high-order PDE operators, we can derive the symbol of the p-multigrid two-grid error propagation operator.

The two-grid multigrid error propagation operator is given by
\begin{equation}
\mathbf{E}_{\text{2mg}} = \mathbf{S}_f \left( \mathbf{I} - \mathbf{P}_{\text{ctof}} \mathbf{A}_c^{-1} \mathbf{R}_{\text{ftoc}} \mathbf{A}_f \right) \mathbf{S}_f,
\end{equation}
where $\mathbf{A}_f$ represents the action of the PDE operator on the fine grid, $\mathbf{S}_f$ represents the fine grid smoother error propagation operator, and $\mathbf{A}_c^{-1}$ represents the coarse grid solve, which may be another multigrid cycle.
$\mathbf{P}_{\text{ctof}}$ and $\mathbf{R}_{\text{ftoc}}$ represent the grid prolongation and restriction operators, respectively.

Some multigrid implementations allow the number of pre and post smoothing passes to be set independently.
This derivation of LFA for p-multigrid allows independently setting these parameters; however, we omit this option in the notation for simplicity.

This error propagation operator can represent both h-multigrid and p-multigrid, depending upon the grid transfer operators and coarse grid representation chosen, but we focus on p-multigrid grid transfer operators in this paper.

\subsubsection{Grid Transfer Operator}\label{sec:grids}

In p-multigrid, grid transfer operators can be represented elementwise and can thus be easily represented in the form of \cref{efficienthighorder}.

\begin{figure}[!ht]
  \centering
  \includegraphics[width=0.48\textwidth]{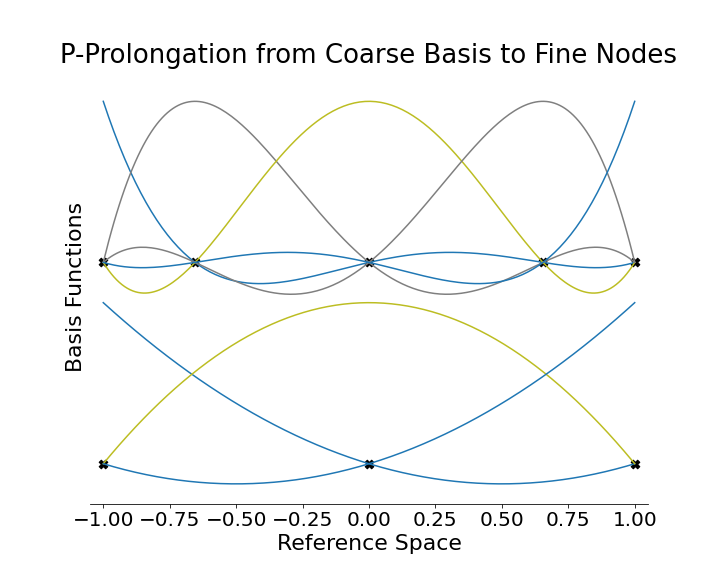}
  \caption{P-Prolongation from Coarse Basis to Fine Basis Points in 1D.}
  \label{fig:p_prolongation}
\end{figure}

The prolongation operator from the coarse to the fine grid interpolates low-order basis functions at the nodes for the high-order basis functions.
\Cref{fig:p_prolongation} shows the evaluation of second order basis function on the Gauss-Lobatto nodes for a fourth order basis.

The prolongation operator from the coarse to the fine grid can be represented by
\begin{equation}
\begin{split}
\mathbf{P}_{\text{ctof}} = \mathbf{P}_f^T \mathbf{P}_e \mathbf{P}_c\\
\mathbf{P}_e = \mathbf{I} \mathbf{D}_{\text{scale}} \mathbf{B}_{\text{ctof}}
\end{split}
\end{equation}
where $\mathbf{B}_{ctof}$ is an interpolation operator from the coarse grid basis to the fine grid basis, $\mathbf{P}_f$ is the fine grid element assembly operator, $\mathbf{P}_c$ is the coarse grid element assembly operator, and $\mathbf{D}_{\text{scale}}$ is a scaling operator to account for node multiplicity across element interfaces.
Restriction from the fine grid to the coarse grid is given by the transpose, $\mathbf{R}_{\text{ftoc}} = \mathbf{P}_{\text{ctof}}^T$.

Following the derivation from \cref{sec:highorder}, we can derive the symbols of $\mathbf{P}_{\text{ctof}}$ and $\mathbf{R}_{\text{ftoc}}$.

\begin{definition}\label{def:prolongation_symbol}
The symbol of the p-prolongation is given by
\begin{equation}
\tilde{\mathbf{P}}_{\text{ctof}} \left( \boldsymbol{\theta} \right) = \mathbf{Q}_f^T \left( \mathbf{P}_e \odot \left[ e^{\imath \left( \mathbf{x}_{j, c} - \mathbf{x}_{i, f} \right) \cdot \boldsymbol{\theta} / \mathbf{h}} \right] \right) \mathbf{Q}_c,
\end{equation}
where $i \in \lbrace 1, 2, \dots, n \cdot \left( p_{\text{fine}} + 1 \right)^d \rbrace$, $\mathbf{h}$ is the length of the element in each dimension, $j \in \lbrace 1, 2, \dots, n \cdot \left( p_{\text{coarse}} + 1 \right)^d \rbrace$, $n$ is the number of components, $p_{\text{fine}}$ and $p_{\text{coarse}}$ are the polynomial degrees of the fine and coarse grid discretizations, respectively, and $d$ is the dimension of the finite element basis.
The matrices $\mathbf{Q}_f$ and $\mathbf{Q}_c$ are the localization mappings for the fine and coarse grid, respectively, and the element p-prolongation operator is given by $\mathbf{P}_e = \mathbf{D}_{\text{scale}} \mathbf{B}_{\text{ctof}}$.
The nodes $\mathbf{x}_{j, c}$ and are $\mathbf{x}_{i, f}$ are on the coarse and fine grids, respectively.
\end{definition}

\begin{definition}\label{def:restriction_symbol}
The symbol of p-restriction is given by the expression
\begin{equation}
\tilde{\mathbf{R}}_{\text{ftoc}} \left( \boldsymbol{\theta} \right) = \mathbf{Q}_c^T \left( \mathbf{R}_e \odot \left[ e^{\imath \left( \mathbf{x}_{j, f} - \mathbf{x}_{i, c} \right) \cdot \boldsymbol{\theta} / \mathbf{h}} \right] \right) \mathbf{Q}_f,
\end{equation}
where $i \in \lbrace 1, 2, \dots, n \cdot \left( p_{\text{coarse}} + 1 \right)^d \rbrace$, $\mathbf{h}$ is the length of the element in each dimension, $j \in \lbrace 1, 2, \dots, n \cdot \left( p_{\text{fine}} + 1 \right)^d \rbrace$, $n$ is the number of components, $p_{\text{fine}}$ and $p_{\text{coarse}}$ are the polynomial degrees of the fine and coarse grid discretizations, respectively, and $d$ is the dimension of the finite element basis.
The matrices $\mathbf{Q}_f$ and $\mathbf{Q}_c$ are the localization mappings for the fine and coarse grid, respectively, and the element p-restriction operator is given by $\mathbf{R}_e = \mathbf{P}_e^T = \mathbf{B}_{\text{ctof}}^T \mathbf{D}_{\text{scale}}$.
The nodes $\mathbf{x}_{j, f}$ and are $\mathbf{x}_{i, c}$ are on the fine and coarse grids, respectively.
\end{definition}

\subsubsection{Multigrid Error Propagation Symbol}\label{sec:multigridsymbol}

Based on the previous discussion, we can combine the expressions in p-multigrid to give the symbol of the error propagation operator.

\begin{definition}\label{def:pmultigrid_symbol}
The symbol for the p-multigrid error propagation operator is given by
\begin{equation}
\tilde{\mathbf{E}}_{\text{2mg}} \left( \nu, \omega, \boldsymbol{\theta} \right) = \tilde{\mathbf{S}}_f \left( \nu, \omega, \boldsymbol{\theta} \right) \left[ \mathbf{I} - \tilde{\mathbf{P}}_{\text{ctof}} \left( \boldsymbol{\theta} \right) \tilde{\mathbf{A}}_c^{-1} \left( \boldsymbol{\theta} \right) \tilde{\mathbf{R}}_{\text{ftoc}} \left( \boldsymbol{\theta} \right) \tilde{\mathbf{A}}_f \left( \boldsymbol{\theta} \right) \right] \tilde{\mathbf{S}}_f \left( \nu, \omega, \boldsymbol{\theta} \right),
\end{equation}
where $\tilde{\mathbf{A}}$ is the symbol of the fine grid operator, $\tilde{\mathbf{S}}_f \left( \nu, \omega, \boldsymbol{\theta} \right)$ is the symbol of the smoother, and $\tilde{\mathbf{A}}_c^{-1} \left( \boldsymbol{\theta} \right)$ is the symbol of the coarse grid solve.
$\tilde{\mathbf{P}}_{\text{ctof}} \left( \boldsymbol{\theta} \right)$ and $\tilde{\mathbf{R}}_{\text{ftoc}} \left( \boldsymbol{\theta} \right)$ represent the symbol of the grid prolongation and restriction operators, respectively.
\end{definition}

Note again that this derivation is applicable for any PDE with a weak form that can be represented by \cref{efficienthighorder}.
This expression can be extended to represent multi-level schemes by recursively applying $\tilde{\mathbf{E}}_{\text{2mg}}$ until the coarsest grid is reached.

\subsection{Extension to H-Multigrid}\label{sec:previouswork}

By using macro-elements, \cite{kumar2019local, brown2019local}, our LFA of p-multigrid can be extended to LFA of h-multigrid.
For 1D analysis, a macro-element is single element comprising of a pair of sub-elements with separate quadrature spaces, which is equivalent to partially assembling the finite element operator on two element subdomains.
Prolongation between a coarse grid element and a fine grid macro-element is given by evaluating the coarse grid basis functions on the fine grid macro-element nodes with the multiplicity correction, as shown in \cref{sec:grids}.
LFA in higher dimensions is given by tensor-products of 1D basis operations, as before.

\begin{figure}[!ht]
  \centering
  \includegraphics[width=0.48\textwidth]{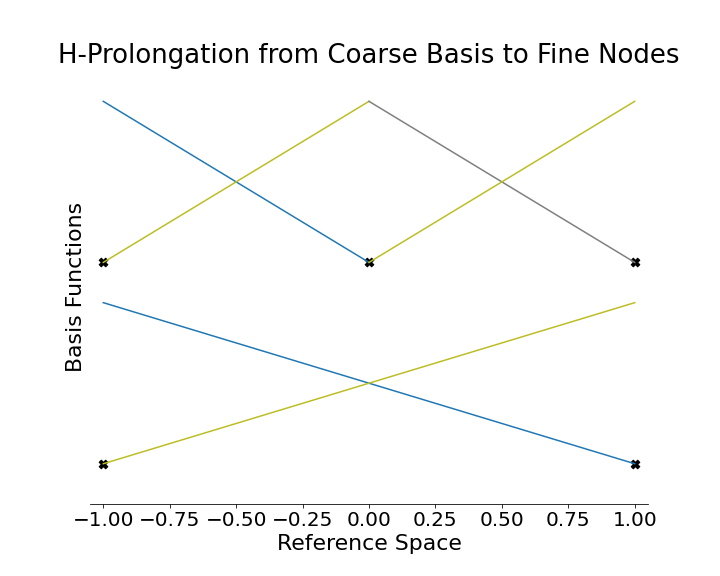}
  \caption{H-Prolongation from Coarse Basis to Fine Basis Points in 1D.}
  \label{fig:h_prolongation}
\end{figure}

In \cref{fig:h_prolongation}, we see an example of interpolation from a coarse grid basis to a fine grid macro-element basis.
The linear shape functions are evaluated at the nodes of the fine grid macro-element, which is a pair of linear sub-elements.
On the left linear micro-element, there are two basis functions which are zero over the domain of the right sub-element, and the reverse is also true.
The element level operator, $\mathbf{A}_e$, on the fine grid macro-element assembles the element level operator for each sub-element into the action of the PDE operator on the full macro-element.

For example, the fine and coarse grid element operators for the one dimensional scalar Laplacian on quadratic elements are given by
\begin{equation}
\mathbf{A}_{e, f} =
\frac{2}{3h}
\begin{bmatrix}
 7  &  -8  &   1  &   0  &   0  \\
-8  &  16  &  -8  &   0  &   0  \\
 1  &  -8  &  14  &  -8  &   1  \\
 0  &   0  &  -8  &  16  &  -8  \\
 0  &   0  &   1  &  -8  &   7  \\
\end{bmatrix},
\end{equation}
where $h$ is the length of the macro-element, and
\begin{equation}
\mathbf{A}_{e, c} =
\frac{1}{3h}
\begin{bmatrix}
 7  &  -8  &   1  \\
-8  &  16  &  -8  \\
 1  &  -8  &   7  \\
\end{bmatrix}.
\end{equation}
The corresponding symbols are computed as in \cref{def:high_order_symbol} above.

Following the derivation from \cref{sec:grids}, we can derive the symbols of grid transfer operators for h-multigrid.

\begin{definition}\label{def:h_prolongation_symbol}
The symbol of the h-prolongation operator is given by
\begin{equation}
\tilde{\mathbf{P}}_{\text{ctof}} \left( \boldsymbol{\theta} \right) = \mathbf{Q}_f^T \left( \mathbf{P}_e \odot \left[ e^{\imath \left( \mathbf{x}_{j, c} - \mathbf{x}_{i, f} \right) \cdot \mathbf{\theta} / \mathbf{h}} \right] \right) \mathbf{Q}_c,
\end{equation}
where $i \in \left\lbrace 1, 2, \dots, n \left( m p + 1 \right)^d \right\rbrace$, $j \in \left\lbrace 1, 2, \dots, n \left( p + 1 \right)^d \right\rbrace$, $\mathbf{h}$ is the length of the macro-element, $d$ is the dimension of the finite element bases, $n$ is the number of components, $m$ is the coarsening factor or the number of micro-elements in each fine grid macro-element, $p$ is the polynomial degree of the discretization, and $d$ is the dimension of the finite element basis.
The matrices $\mathbf{Q}_f$ and $\mathbf{Q}_c$ are the localization mappings for the fine and coarse grid, respectively, and the macro-element h-prolongation operator is given by $\mathbf{P}_e = \mathbf{I} \mathbf{D}_{\text{scale}} \mathbf{B}_{\text{ctof}}$.
\end{definition}

\begin{definition}\label{def:h_restriction_symbol}
The symbol of h-restriction operator is given by
\begin{equation}
\tilde{\mathbf{R}}_{\text{ftoc}} \left( \boldsymbol{\theta} \right) = \mathbf{Q}_c^T \left( \mathbf{R}_e \odot \left[ e^{\imath \left( \mathbf{x}_{j, f} - \mathbf{x}_{i, c} \right) \cdot \boldsymbol{\theta} / \mathbf{h}} \right] \right) \mathbf{Q}_f,
\end{equation}
where $i \in \left\lbrace 1, 2, \dots, n \left( p + 1 \right)^d \right\rbrace$, $j \in \left\lbrace 1, 2, \dots, n \left( m p + 1 \right)^d \right\rbrace$, $\mathbf{h}$ is the length of the macro-element, $d$ is the dimension of the finite element bases, $n$ is the number of components, $m$ is the coarsening factor or the number of micro-elements in each fine grid macro-element, $p$ is the polynomial degree of the discretization, and $d$ is the dimension of the finite element basis.
The matrices $\mathbf{Q}_f$ and $\mathbf{Q}_c$ are the localization mappings for the fine and coarse grid, respectively, and the macro-element h-restriction operator is given by $\mathbf{R}_e = \mathbf{P}_e^T = \mathbf{B}_{\text{ctof}}^T \mathbf{D}_{\text{scale}} \mathbf{I}$.
\end{definition}

By representing the fine grid with macro-elements and the prolongation operator with this interpolation, this LFA of p-multigrid exactly reproduces the results of He and MacLachlan \cite{he2020two} for LFA of high-order h-multigrid.

Our LFA of h-multigrid presented here is based on a specific basis of the Fourier space used in \cite{kumar2019local} rather than the commonly used Fourier modes, see \cite{MR1807961, wienands2004practical}, where the symbol of each component in multigrid methods is formed based on different harmonic frequencies.
The symbols of grid-transfer operators described in \cref{def:h_prolongation_symbol,def:h_restriction_symbol} give a general way for multigrid coarsening with factor $m$, and this framework is simpler, especially for the high-order discretizations described in \cite{he2020two}.
Also, our LFA for h-multigrid is not restricted to finite element discretizations with uniformly spaced nodes.

On uniform rectangular meshes, linear finite elements produce the same discetized operator as finite differencing.
The nine-point stencil for the Laplace operator in two dimensions is given by

\begin{equation}
\frac{1}{3}
\begin{bmatrix}
-1  &  -1  &  -1  \\
-1  &   8  &  -1  \\
-1  &  -1  &  -1  \\
\end{bmatrix}
\end{equation}
with a corresponding local Fourier analysis symbol given by

\begin{equation}
\tilde{A} \left( \theta_1, \theta_2 \right) = \frac{2}{3} \left( 4 - \cos \left( \theta_1 \right) - \cos \left( \theta_2 \right) - 2 \cos \left( \theta_1 \right) \cos \left( \theta_2 \right) \right).
\end{equation}

The assembled matrix for a single linear element in two dimensions is given by

\begin{equation}
\mathbf{A}_e =
\frac{1}{3}
\begin{bmatrix}
 2    &  -1/2  &  -1/2  &  -1    \\
-1/2  &   2    &  -1    &  -1/2  \\
-1/2  &  -1    &   2    &  -1/2  \\
-1    &  -1/2  &  -1/2  &   2    \\
\end{bmatrix},
\end{equation}
with a corresponding local Fourier analysis symbol given by

\begin{equation}
\begin{split}
\tilde{A} \left( \theta_1, \theta_2 \right) & = \mathbf{Q}^T \left( \mathbf{A}_e \odot \left[ e^{\imath \left( \mathbf{x}_j - \mathbf{x}_i \right) \cdot \boldsymbol{\theta}} \right] \right) \mathbf{Q}\\ & = \frac{2}{3} \left( 4 - \cos \left( \theta_1 \right) - \cos \left( \theta_2 \right) - 2 \cos \left( \theta_1 \right) \cos \left( \theta_2 \right) \right),
\end{split}
\end{equation}
where $\mathbf{Q} = \begin{bmatrix} 1 & 1 & 1 & 1 \end{bmatrix}^T$.

We can use the LFA of multigrid given by \cref{def:pmultigrid_symbol} with the prolongation and restriction symbols given by \cref{def:h_prolongation_symbol,def:h_restriction_symbol} to reproduce LFA of h-multigrid methods for finite differencing where the stencil can be represented by a finite element discretization.
This LFA of arbitrary second-order PDEs with high-order finite element discretizations agrees with previous work on LFA of PDE operators derived with finite differencing with analogous stencils.

As mentioned before, our focus of this work is p-multigrid, so we will not expand the discussion of h-multigrid further method here.
Applying this LFA framework to h-multigrid or hp-multigrid methods are topics for future research.

\section{Numerical Results}\label{sec:results}

In this section, we present numerical results for this analysis for the scalar Laplacian in one and two dimensions with $H^1$ Lagrange bases on Gauss-Lobatto points with Gauss-Legendre quadrature.
Next, we validate these results with numerical experiments for the scalar Laplacian in three dimensions.
Lastly, we consider linear elasticity in three dimensions.

\subsection{Scalar Laplacian - 1D Convergence Factors}\label{sec:1dresults}

\subsubsection{Jacobi Smoothing}

\begin{figure}[!tbp]
  \centering
    \subfloat[Convergence for $p = 4$ to $p = 2$, $\nu = 1$]{\includegraphics[width=0.48\textwidth]{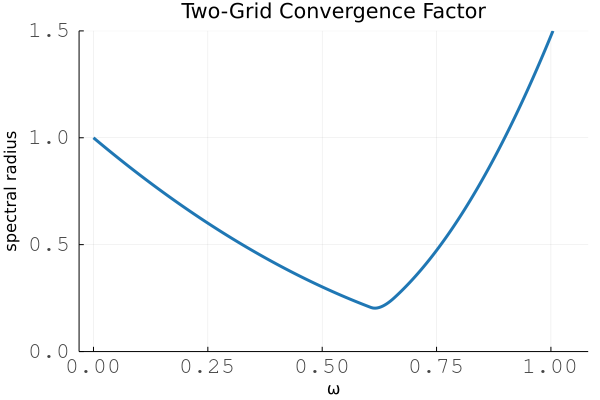}\label{fig:two_grid_5_3}}
    \subfloat[Convergence for $p = 4$ to $p = 1$, $\nu = 1$]{\includegraphics[width=0.48\textwidth]{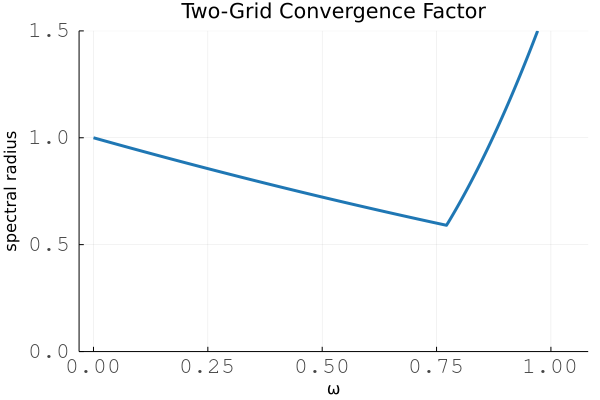}\label{fig:two_grid_5_2}} \\
    \subfloat[Convergence for $p = 4$ to $p = 2$, $\nu = 2$]{\includegraphics[width=0.48\textwidth]{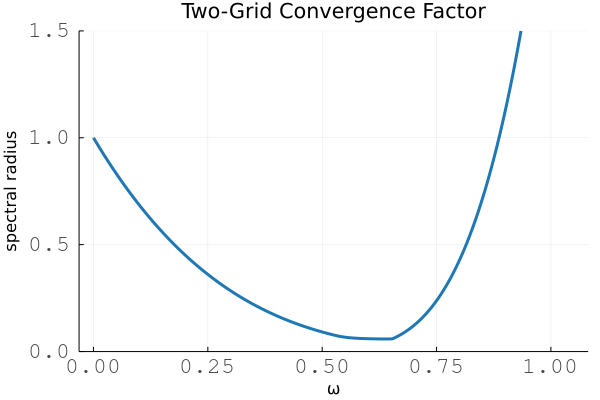}\label{fig:two_grid_5_3_2smooth}}
    \subfloat[Convergence for $p = 4$ to $p = 1$, $\nu = 2$]{\includegraphics[width=0.48\textwidth]{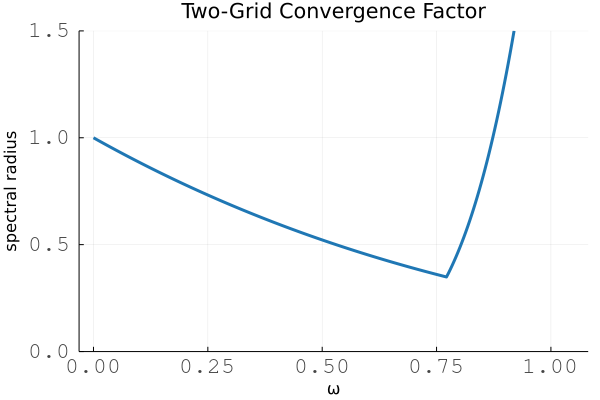}\label{fig:two_grid_5_2_2smooth}} \\
  \caption{Two-grid analysis for Jacobi smoothing for high-order finite elements for the 1D Laplacian.}
\end{figure}

In \cref{fig:two_grid_5_3} and \cref{fig:two_grid_5_2}, we plot the two-grid convergence factor for p-multigrid with a single iteration of Jacobi pre and post-smoothing for the one dimensional Laplacian as a function of the Jacobi smoothing parameter $\omega$,
and in \cref{fig:two_grid_5_3_2smooth} and \cref{fig:two_grid_5_2_2smooth} we plot the two-grid convergence factor for p-multigrid with two iterations of Jacobi pre and post-smoothing for the one dimensional Laplacian as a function of the Jacobi smoothing parameter $\omega$.
On the left we show conservative coarsening from quartic to quadratic elements and on the right we show more aggressive coarsening from quartic to linear elements.
As expected, the two-grid convergence factor decreases as we coarsen more rapidly.
Also, the effect of underestimating the optimal Jacobi smoothing parameter, $\omega$, is less pronounced than the effect of overestimating the smoothing parameter, especially with a higher number of pre and post-smooths.

In contrast to the previous work on h-multigrid for high-order finite elements, \cite{he2020two}, poorly chosen values of $\omega < 1.0$ can result in a spectral radius of the p-multigrid error propagation symbol that is greater than $1$, indicating that application of p-multigrid with Jacobi smoothing at these parameter values will result in increased error.

\begin{table}[ht!]
\begin{center}
\begin{tabular}{l cc cc cc}
  \toprule
  $p_{\text{fine}}$ to $p_{\text{coarse}}$  &  \multicolumn{2}{c}{$\nu = 1$}  &  \multicolumn{2}{c}{$\nu = 2$}  &  \multicolumn{2}{c}{$\nu = 3$}  \\
                       &  $\rho_{\min}$ & $\omega_{\text{opt}}$  &  $\rho_{\min}$ & $\omega_{\text{opt}}$  &  $\rho_{\min}$ & $\omega_{\text{opt}}$  \\
  \toprule
  $p = 2$ to $p = 1$   &  0.137 & 0.63  &  0.060 & 0.69  &  0.041 & 0.72   \\
  \midrule
  $p = 4$ to $p = 2$   &  0.204 & 0.62  &  0.059 & 0.64  &  0.045 & 0.70   \\
  $p = 4$ to $p = 1$   &  0.591 & 0.77  &  0.350 & 0.77  &  0.207 & 0.77   \\
  \midrule
  $p = 8$ to $p = 4$   &  0.250 & 0.60  &  0.068 & 0.60  &  0.033 & 0.63   \\
  $p = 8$ to $p = 2$   &  0.668 & 0.73  &  0.446 & 0.73  &  0.298 & 0.73   \\
  $p = 8$ to $p = 1$   &  0.874 & 0.78  &  0.764 & 0.78  &  0.668 & 0.78   \\
  \midrule
  $p = 16$ to $p = 8$  &  0.300 & 0.57  &  0.090 & 0.57  &  0.035 & 0.58   \\
  $p = 16$ to $p = 4$  &  0.719 & 0.69  &  0.517 & 0.69  &  0.371 & 0.69   \\
  $p = 16$ to $p = 2$  &  0.906 & 0.73  &  0.820 & 0.73  &  0.743 & 0.73   \\
  $p = 16$ to $p = 1$  &  0.968 & 0.74  &  0.936 & 0.74  &  0.906 & 0.74   \\
  \bottomrule
\end{tabular}
\end{center}
\caption{Two-grid convergence factor and optimal Jacobi parameter for the 1D Laplacian.}
\label{table:two_grid_1d}
\end{table}

The results in \cref{table:two_grid_1d} provide the LFA convergence factor and optimal values of $\omega$ for two-grid high-order p-multigrid for a variety of basis polynomial degrees and coarsening factors.

Traditional estimates of the optimal smoothing parameter based upon extremal eigenvalues of the preconditioned operator are incompatible with this LFA framework.
Optimal parameter estimation is an open question for high-order p-multigrid, but optimization techniques, such as those discussed in \cite{brown2021tuning}, can be used to tune these parameters, especially for more complex PDEs.

\subsubsection{Chebyshev Smoothing}

\begin{figure}[!tbp]
  \centering
    \subfloat[Convergence for $p = 4$ to $p = 2$, $\nu = 1$]{\includegraphics[width=0.48\textwidth]{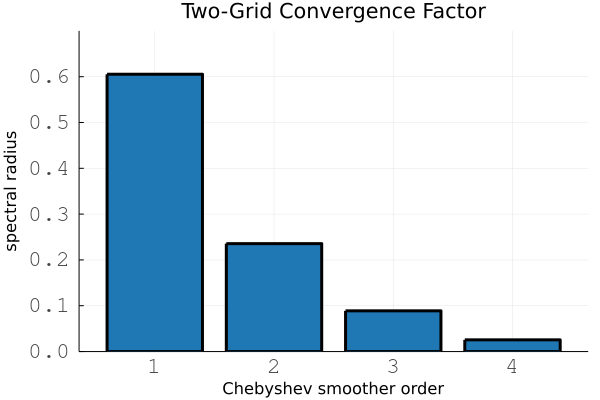}\label{fig:two_grid_converge_4_to_2_chebyshev}}
    \subfloat[Convergence for $p = 4$ to $p = 1$, $\nu = 1$]{\includegraphics[width=0.48\textwidth]{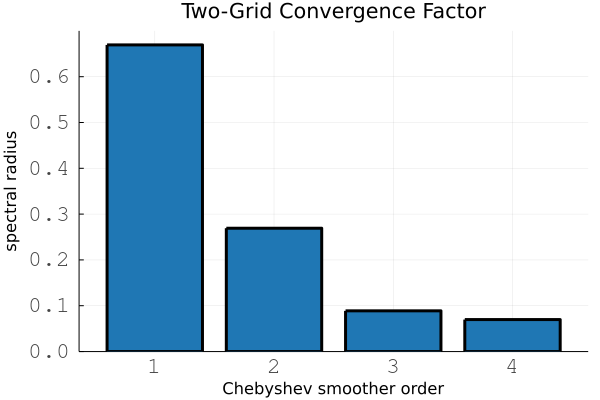}\label{fig:two_grid_converge_4_to_1_chebyshev}}
  \caption{Two-grid analysis for Chebyshev smoothing for high-order finite elements for the 1D Laplacian.}
\end{figure}

In \cref{fig:two_grid_converge_4_to_2_chebyshev,fig:two_grid_converge_4_to_1_chebyshev} we plot the two-grid convergence factor for p-multigrid with Chebyshev pre and post-smoothing for the one dimensional Laplacian as a function of the Chebyshev order, $k$.
On the left we show conservative coarsening from quartic to quadratic elements and on the right we show more aggressive coarsening from quartic to linear elements.
As expected, the two-grid convergence factor degrades as we coarsen more rapidly.

\begin{table}[ht!]
\begin{center}
\begin{tabular}{l c c c c}
  \toprule
  $p_{\text{fine}}$ to $p_{\text{coarse}}$  &  $k = 1$   &  $k = 2$   &  $k = 3$   &  $k = 4$   \\
  \toprule
  $p = 2$ to $p = 1$   &  0.545  &  0.220  &  0.063  &  0.017  \\
  \midrule
  $p = 4$ to $p = 2$   &  0.576  &  0.222  &  0.089  &  0.025  \\
  $p = 4$ to $p = 1$   &  0.623  &  0.269  &  0.089  &  0.070  \\
  \midrule
  $p = 8$ to $p = 4$   &  0.638  &  0.244  &  0.074  &  0.022  \\
  $p = 8$ to $p = 2$   &  0.657  &  0.260  &  0.097  &  0.059  \\
  $p = 8$ to $p = 1$   &  0.881  &  0.674  &  0.510  &  0.393  \\
  \midrule
  $p = 16$ to $p = 8$  &  0.664  &  0.253  &  0.075  &  0.022  \\
  $p = 16$ to $p = 4$  &  0.714  &  0.328  &  0.135  &  0.059  \\
  $p = 16$ to $p = 2$  &  0.907  &  0.741  &  0.602  &  0.496  \\
  $p = 16$ to $p = 1$  &  0.970  &  0.912  &  0.857  &  0.809  \\
  \bottomrule
\end{tabular}
\end{center}
\caption{Two-grid convergence factor with Chebyshev smoothing for 1D Laplacian.}
\label{table:two_grid_1d_chebyshev}
\end{table}

The results in \cref{table:two_grid_1d_chebyshev} provide the LFA-predicted convergence factor and optimal values of $k$ for two-grid high-order p-multigrid for a variety of coarsening rates and orders of Chebyshev smoother.
From this table, we can see that the effectiveness of higher order Chebyshev smoothers degrades as we coarsen more aggressively, but Chebyshev smoothing still provides better two-grid convergence than multiple pre and post-smoothing Jacobi iterations.

\begin{table}[ht!]
\begin{center}
\begin{tabular}{l c c c c}
  \toprule
  \multicolumn{5}{c}{$\lambda_{\min} = 0.2 \hat{\lambda}_{\max}$} \\
  \toprule
  $p_{\text{fine}}$ to $p_{\text{coarse}}$  &  $k = 1$   &  $k = 2$   &  $k = 3$   &  $k = 4$   \\
  \toprule
  $p = 4$ to $p = 2$   &  0.410  &  0.093  &  0.043  &  0.024  \\
  $p = 4$ to $p = 1$   &  0.611  &  0.250  &  0.106  &  0.071  \\
  \midrule
  $p = 8$ to $p = 4$   &  0.435  &  0.081  &  0.016  &  0.007  \\
  $p = 8$ to $p = 1$   &  0.891  &  0.739  &  0.623  &  0.529  \\
  \midrule
  $p = 16$ to $p = 8$  &  0.443  &  0.081  &  0.015  &  0.006  \\
  $p = 16$ to $p = 1$  &  0.973  &  0.931  &  0.894  &  0.861  \\
  \toprule
  \multicolumn{5}{c}{$\lambda_{\min} = 0.3 \hat{\lambda}_{\max}$} \\
  \toprule
  $p_{\text{fine}}$ to $p_{\text{coarse}}$  &  $k = 1$   &  $k = 2$   &  $k = 3$   &  $k = 4$   \\
  \toprule
  $p = 4$ to $p = 2$   &  0.279  &  0.070  &  0.042  &  0.031  \\
  $p = 4$ to $p = 1$   &  0.638  &  0.332  &  0.184  &  0.104  \\
  \midrule
  $p = 8$ to $p = 4$   &  0.289  &  0.050  &  0.023  &  0.012  \\
  $p = 8$ to $p = 1$   &  0.899  &  0.777  &  0.682  &  0.599  \\
  \midrule
  $p = 16$ to $p = 8$  &  0.294  &  0.055  &  0.020  &  0.010  \\
  $p = 16$ to $p = 1$  &  0.975  &  0.942  &  0.913  &  0.885  \\
  \bottomrule
\end{tabular}
\end{center}
\caption{Two-grid convergence factor with Chebyshev smoothing for 1D Laplacian with modified lower eigenvalue bound.}
\label{table:two_grid_1d_chebyshev_eigenvalues}
\end{table}

The results in \cref{table:two_grid_1d_chebyshev_eigenvalues} provide the LFA-predicted convergence factor for two-grid high-order p-multigrid for a variety of coarsening rates and orders of Chebyshev smoother with different scaling factors for the minimum eigenvalue estimate used in the Chebyshev iterations.
Increasing the lower eigenvalue estimate results in the Chebyshev method better targeting high frequency error modes, which results in improved two-grid convergence when halving the polynomial degree of the basis functions.
However, increasing the lower eigenvalue estimate results in worse two-grid convergence for aggressive coarsening directly to linear elements.

\subsection{Scalar Laplacian - 2D Convergence Factors}\label{sec:2dresults}

\subsubsection{Jacobi Smoothing}

\begin{figure}[!tbp]
  \centering
  \subfloat[Convergence for $p = 4$ to $p = 1$, $\nu = 1$]{\includegraphics[width=0.48\textwidth]{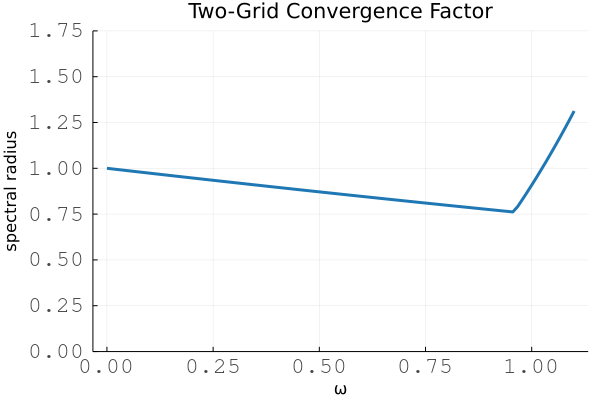}\label{fig:two_grid_converge_4_to_1_1_smooth_2d}}
  \hfill
  \subfloat[Convergence for $p = 4$ to $p = 1$, $\nu = 2$]{\includegraphics[width=0.48\textwidth]{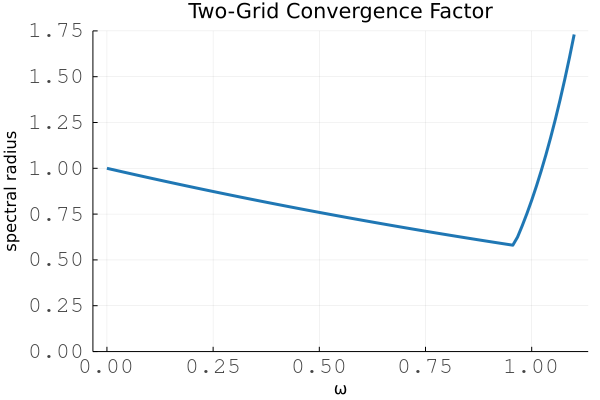}\label{fig:two_grid_converge_4_to_1_2_smooth_2d}}
  \caption{Convergence for high-order finite elements for the 2D Laplacian.}
\end{figure}

In \Cref{fig:two_grid_converge_4_to_1_1_smooth_2d} and \cref{fig:two_grid_converge_4_to_1_2_smooth_2d} we show the two-grid convergence factor for p-multigrid with one and two iterations of Jacobi smoothing, respectively, for the two dimensional Laplacian as a function of the Jacobi smoothing parameter $\omega$.

\begin{table}[ht!]
\begin{center}
\begin{tabular}{l cc cc cc}
  \toprule
  $p_{\text{fine}}$ to $p_{\text{coarse}}$  &  \multicolumn{2}{c}{$\nu = 1$}  &  \multicolumn{2}{c}{$\nu = 2$}  &  \multicolumn{2}{c}{$\nu = 3$}  \\
                      &  $\rho_{\min}$  &  $\omega_{\text{opt}}$  &  $\rho_{\min}$ & $\omega_{\text{opt}}$  &  $\rho_{\min}$ & $\omega_{\text{opt}}$  \\
  \toprule
  $p = 2$ to $p = 1$  &  0.230 & 0.95  &  0.091 & 0.99  &  0.061 & 1.03   \\
  \midrule
  $p = 4$ to $p = 2$  &  0.388 & 0.82  &  0.151 & 0.82  &  0.078 & 0.83   \\
  $p = 4$ to $p = 1$  &  0.763 & 0.95  &  0.582 & 0.95  &  0.444 & 0.95   \\
  \midrule
  $p = 8$ to $p = 4$  &  0.646 & 0.79  &  0.418 & 0.79  &  0.272 & 0.79   \\
  $p = 8$ to $p = 2$  &  0.858 & 0.84  &  0.737 & 0.84  &  0.633 & 0.84   \\
  $p = 8$ to $p = 1$  &  0.952 & 0.87  &  0.907 & 0.87  &  0.864 & 0.87   \\
  \bottomrule
\end{tabular}
\end{center}
\caption{Two-grid convergence factor and optimal Jacobi parameter for 2D Laplacian.}
\label{table:two_grid_2d}
\end{table}

The results in \cref{table:two_grid_2d} provide the LFA-predicted convergence factor and optimal values of $\omega$ for two-grid high-order p-multigrid for a variety of basis polynomial degrees and coarsening factors.

\subsubsection{Chebyshev Smoothing}

\begin{table}[ht!]
\begin{center}
\begin{tabular}{l c c c c}
  \toprule
  $p_{\text{fine}}$ to $p_{\text{coarse}}$  &  $k = 1$   &  $k = 2$   &  $k = 3$   &  $k = 4$   \\
  \toprule
  $p = 2$ to $p = 1$   &  0.621  &  0.252  &  0.075  &  0.039  \\
  \midrule
  $p = 4$ to $p = 2$   &  0.607  &  0.281  &  0.085  &  0.047  \\
  $p = 4$ to $p = 1$   &  0.768  &  0.424  &  0.219  &  0.127  \\
  \midrule
  $p = 8$ to $p = 4$   &  0.669  &  0.278  &  0.110  &  0.055  \\
  $p = 8$ to $p = 2$   &  0.864  &  0.633  &  0.456  &  0.336  \\
  $p = 8$ to $p = 1$   &  0.956  &  0.873  &  0.795  &  0.730  \\
  \midrule
  $p = 16$ to $p = 8$  &  0.855  &  0.613  &  0.435  &  0.319  \\
  $p = 16$ to $p = 4$  &  0.938  &  0.822  &  0.719  &  0.634  \\
  $p = 16$ to $p = 2$  &  0.976  &  0.928  &  0.882  &  0.842  \\
  $p = 16$ to $p = 1$  &  0.992  &  0.975  &  0.959  &  0.944  \\
  \bottomrule
\end{tabular}
\end{center}
\caption{Two-grid convergence factor with Chebyshev smoothing for 2D Laplacian.}
\label{table:two_grid_2d_chebyshev}
\end{table}

The results in \cref{table:two_grid_2d_chebyshev} provide the LFA-predicted convergence factor for two-grid high-order p-multigrid for a variety of coarsening rates and orders of Chebyshev smoother.
The two-grid convergence factor still degrades and the effectiveness of higher order Chebyshev smoothers is again reduced as we coarsen more aggressively.

\begin{table}[ht!]
\begin{center}
\begin{tabular}{l c c c c}
  \toprule
  \multicolumn{5}{c}{$\lambda_{\min} = 0.2 \hat{\lambda}_{\max}$} \\
  \toprule
  $p_{\text{fine}}$ to $p_{\text{coarse}}$  &  $k = 1$   &  $k = 2$   &  $k = 3$   &  $k = 4$   \\
  \toprule
  $p = 4$ to $p = 2$   &  0.450  &  0.137  &  0.067  &  0.050  \\
  $p = 4$ to $p = 1$   &  0.786  &  0.525  &  0.362  &  0.255  \\
  \midrule
  $p = 8$ to $p = 4$   &  0.668  &  0.330  &  0.172  &  0.106  \\
  $p = 8$ to $p = 1$   &  0.960  &  0.899  &  0.848  &  0.801  \\
  \toprule
  \multicolumn{5}{c}{$\lambda_{\min} = 0.3 \hat{\lambda}_{\max}$} \\
  \toprule
  $p_{\text{fine}}$ to $p_{\text{coarse}}$  &  $k = 1$   &  $k = 2$   &  $k = 3$   &  $k = 4$   \\
  \toprule
  $p = 4$ to $p = 2$   &  0.407  &  0.106  &  0.073  &  0.059  \\
  $p = 4$ to $p = 1$   &  0.803  &  0.590  &  0.447  &  0.341  \\
  \midrule
  $p = 8$ to $p = 4$   &  0.691  &  0.409  &  0.256  &  0.164  \\
  $p = 8$ to $p = 1$   &  0.963  &  0.915  &  0.874  &  0.835  \\
  \bottomrule
\end{tabular}
\end{center}
\caption{Two-grid convergence factor with Chebyshev smoothing for 2D Laplacian with modified lower eigenvalue bound.}
\label{table:two_grid_2d_chebyshev_eigenvalues}
\end{table}

The results in \cref{table:two_grid_2d_chebyshev_eigenvalues} provide the LFA-predicted convergence factor and optimal values of $k$ for two-grid high-order p-multigrid for a variety of coarsening rates and orders of Chebyshev smoother with different scaling factors for the minimum eigenvalue estimate used in the Chebyshev iterations.
As in one dimension, increasing the lower eigenvalue estimate results in better two-grid convergence when conventional coarsening by halving the polynomial degree of the basis functions but results in worse two-grid convergence with aggressive coarsening.

\subsection{Scalar Laplacian - 3D Convergence Factors}\label{sec:3dresults}

In this section, we compare the LFA two-grid convergence factors to numerical results.
Our numerical experiments were conducted using the libCEED \cite{libceed-joss} with PETSc \cite{petsc-user-ref} multigrid example found in the libCEED repository.
PETSc provides the mesh management, linear solvers, and multigrid preconditioner while libCEED provides the matrix-free operator evaluation.

We recover the manufactured solution given by
\begin{equation}
f \left( x, y, z \right) = x y z \sin \left( \pi x \right) \sin \left( \pi \left( 1.23 + 0.5 y \right) \right) \sin \left( \pi \left( 2.34 + 0.25 z \right) \right)
\end{equation}
on the domain $\left[ -3, 3 \right]^3$ with Dirichlet boundary conditions for finite element discretizations with varying orders on approximately 8 million degrees of freedom for a variety of test cases.

Although LFA is defined on infinite grids, which most naturally translate to periodic problems, LFA-predicted convergence factors are also accurate for appropriate problems with other boundary conditions on finite grids after appropriate reformulation of the analysis \cite{rodrigo2019validity}.
Since it is not feasable to tune multigrid individually for each set of boundary conditions in engineering domains, we explore the sharpness and robustness of infinite-grid LFA by comparing directly to the test problem with Dirichlet boundary conditions.

\subsubsection{Jacobi Smoothing}

Since the Chebyshev smoothing is based upon the Jacobi preconditioned operator, it is important to validate the LFA of the Jacobi smoothing before considering Chebyshev smoothing.
We use simple Jacobi smoothing with a weight of $\omega = 1.0$ to validate the LFA.

\begin{table}[ht!]
\begin{center}
\begin{tabular}{l c c}
  \toprule
  $p_{\text{fine}}$ to $p_{\text{coarse}}$  &  LFA  &  libCEED  \\
  \toprule
  $p = 2$ to $p = 1$   &  0.312  &  0.301  \\
  \midrule
  $p = 4$ to $p = 2$   &  1.436  &  1.402  \\
  $p = 4$ to $p = 1$   &  1.436  &  1.401  \\
  \midrule
  $p = 8$ to $p = 4$   &  1.989  &  1.885  \\
  $p = 8$ to $p = 2$   &  1.989  &  1.874  \\
  $p = 8$ to $p = 1$   &  1.989  &  1.875  \\
  \bottomrule
\end{tabular}
\end{center}
\caption{LFA and experimental two-grid convergence factors with Jacobi smoothing for 3D Laplacian with $\omega = 1.0$.}
\label{table:two_grid_3d_jacobi}
\end{table}

The results in \cref{table:two_grid_3d_jacobi} provide the LFA and experimental convergence factors for the test problem.
As expected, the high-order fine grid problems diverge with a smoothing factor of $\omega = 1.0$; however, the LFA provides reasonable upper bounds on the true convergence factor seen in the experimental results.

\subsubsection{Chebyshev Smoothing}

We used the LFA estimates of the maximal eigenvalue to set the extremal eigenvalues used the Chebyshev iteration in PETSc, using $\lambda_{\text{min}} = 0.1 \hat{\lambda}_{\text{max}}$ and $\lambda_{\text{max}} = 1.0 \hat{\lambda}_{\text{max}}$, where $\hat{\lambda}_{\text{max}}$ is the estimated maximal eigenvalue of the symbol of the Jacobi preconditioned operator.

\begin{table}[ht!]
\begin{center}
\begin{tabular}{l cc cc cc}
  \toprule
  $p_{\text{fine}}$ to $p_{\text{coarse}}$  &  \multicolumn{2}{c}{$k = 2$}  &  \multicolumn{2}{c}{$k = 3$}  &  \multicolumn{2}{c}{$k = 4$}  \\
                      &  LFA  &  libCEED  &  LFA  &  libCEED  &  LFA  &  libCEED  \\
  \toprule
  $p = 2$ to $p = 1$  &  0.253 & 0.222  &  0.076 & 0.058  &  0.041 & 0.033  \\
  \midrule
  $p = 4$ to $p = 2$  &  0.277 & 0.251  &  0.111 & 0.097  &  0.062 & 0.050  \\
  $p = 4$ to $p = 1$  &  0.601 & 0.587  &  0.416 & 0.398  &  0.295 & 0.276  \\
  \midrule
  $p = 8$ to $p = 4$  &  0.398 & 0.391  &  0.197 & 0.195  &  0.121 & 0.110  \\
  $p = 8$ to $p = 2$  &  0.748 & 0.743  &  0.611 & 0.603  &  0.506 & 0.469  \\
  $p = 8$ to $p = 1$  &  0.920 & 0.914  &  0.871 & 0.861  &  0.827 & 0.814  \\
  \bottomrule
\end{tabular}
\end{center}
\caption{LFA and experimental two-grid convergence factors with Chebyshev smoothing for 3D Laplacian.}
\label{table:two_grid_3d_chebyshev}
\end{table}

\cref{table:two_grid_3d_chebyshev} shows that the LFA predictions agree well with the experimental convergence factors.
As with the one and two dimensional results, rapid coarsening of the polynomial degree of the bases decreases the effectiveness of higher order Chebyshev smoothing.

\subsection{Linear Elasticity - 3D Convergence Factors}\label{sec:solidsresults}

To demonstrate the suitability of this LFA formulation for more complex PDEs, we consider linear elasticity in three dimensions.
The strong form of the static balance of linear momentum at small strain for the three dimensional linear elasticity problem is given by \cite{hughes2012finite} as
\begin{equation}
\nabla \cdot \boldsymbol{\sigma} + \boldsymbol{g} = \boldsymbol{0},
\end{equation}
where $\boldsymbol{\sigma}$ is the stress function and $\boldsymbol{g}$ is the forcing function.
This strong form has the corresponding weak form
\begin{equation}
\int_{\Omega} \nabla \mathbf{v} : \boldsymbol{\sigma} dV - \int_{\partial \Omega} \mathbf{v} \cdot \left( \boldsymbol{\sigma} \cdot \hat{\mathbf{n}} \right) dS - \int_{\Omega} \mathbf{v} \cdot \mathbf{g} dV = 0, \forall \mathbf{v} \in \mathcal{V}
\end{equation}
for some displacement $\mathbf{u} \in \mathcal{V} \subset H^1 \left( \Omega \right)$, where $:$ denotes contraction over both components and dimensions.

Linear elasticity constitutive modeling is based upon the Lamé parameters,
\begin{equation}
\lambda = \frac{E \nu}{\left( 1 + \nu \right) \left( 1 - 2 \nu \right)},\quad \mu = \frac{E}{2 \left( 1 + \nu \right)},
\end{equation}
where $E$ is the Young's modulus and $\nu$ is the Poisson's ratio for the materiel.

In the linear elasticity constitutive model, the symmetric strain tensor is given by
\begin{equation}
\boldsymbol{\epsilon} = \frac{1}{2} \left( \nabla \mathbf{u} + \nabla \mathbf{u}^T \right),
\end{equation}
and the linear elasticity constitutive law is given by $\boldsymbol{\sigma} = \mathsf{C} : \boldsymbol{\epsilon}$ where
\begin{equation}
\mathsf{C} =
\begin{bmatrix}
   \lambda + 2\mu & \lambda & \lambda & & & \\
   \lambda & \lambda + 2\mu & \lambda & & & \\
   \lambda & \lambda & \lambda + 2\mu & & & \\
   & & & \mu & & \\
   & & & & \mu & \\
   & & & & & \mu
\end{bmatrix}.
\end{equation}

We can represent this PDE in the form given by \cref{efficienthighorder} and therefore investigate the LFA of p-multgrid for this PDE.

\begin{table}[ht!]
\begin{center}
\begin{tabular}{l c c c c}
  \toprule
  $p_{\text{fine}}$ to $p_{\text{coarse}}$  &  $k = 1$   &  $k = 2$   &  $k = 3$   &  $k = 4$   \\
  \toprule
  $p = 2$ to $p = 1$   &  0.815  &  0.521  &  0.321  &  0.204  \\
  \midrule
  $p = 4$ to $p = 2$   &  0.878  &  0.666  &  0.499  &  0.382  \\
  $p = 4$ to $p = 1$   &  0.969  &  0.908  &  0.850  &  0.800  \\
  \midrule
  $p = 8$ to $p = 4$   &  0.939  &  0.826  &  0.725  &  0.643  \\
  $p = 8$ to $p = 2$   &  0.981  &  0.943  &  0.906  &  0.873  \\
  $p = 8$ to $p = 1$   &  0.995  &  0.984  &  0.973  &  0.963  \\
  \bottomrule
\end{tabular}
\end{center}
\caption{Two-grid convergence factor with Chebyshev smoothing for 3D linear elasticity.}
\label{table:two_grid_3d_linear_elasticity}
\end{table}

The results in \cref{table:two_grid_3d_linear_elasticity} provide the LFA convergence factor and optimal values of $k$ for two-grid high-order p-multigrid for a variety of coarsening rates and orders of Chebyshev smoother.
The LFA two-grid convergence factor for three dimensional linear elasticity is much poorer than the LFA two-grid convergence factor for the three dimensional scalar Laplacian, and the effect of rapid coarsening in the polynomial degree of the basis functions is far more pronounced for linear elasticity.

\section{Conclusions}\label{sec:conclusion}

In this paper we developed LFA of p-multigrid with arbitrary second-order PDEs using high-order finite element discretizations by using an operator representation for efficient application of matrix-free implementations.
We introduced LFAToolkit.jl \cite{thompson2021toolkit}, a new Julia package for LFA of high-order finite element methods, and used LFAToolkit.jl to investigate p-multigrid method with Jacobi and Chebyshev semi-iterative method smoothing and aggressive coarsening strategies.
We observed that the performance of p-multigrid with these two polynomial smoothers degrades as we coarsen more aggressively.

The LFA of p-multigrid framework presented here can be extended to the LFA of h-multigrid methods and reproduces previous work in this area.
Also, we briefly described how this work can be extended to LFA of h-multigrid for finite difference discretization stencils that can be represented as finite difference methods.

The LFA formulation in LFAToolkit.jl accurately predicts the performance of two-grid p-multigrid with polynomial smoothers for the scalar Laplacian in three dimensions and can provide two-grid convergence factors for more complex PDEs, such as linear elasticity in three dimensions.

\bibliographystyle{siamplain}
\bibliography{references}

\end{document}